\documentclass{article}
\usepackage{latexsym}
\usepackage{amsfonts}
\usepackage{amscd}
\usepackage{amsmath, amsthm, amssymb}
\usepackage[all]{xy}
\newtheorem{thm}{Theorem}[section]
\newtheorem{cor}[thm]{Corollary}
\newtheorem{lem}[thm]{Lemma}
\newtheorem{prop}[thm]{Proposition}
\newtheorem{define}{Definition}
\newtheorem{rmk}{Remark}
\newtheorem{example}{Example}

\newcommand{\lie}[1]{\mathfrak{#1}}
\newcommand{\clie}[1]{\mathfrak{#1}_{\C}}
\newcommand{\Z}{\mathbb{Z}}
\newcommand{\R}{\mathbb{R}}
\newcommand{\C}{\mathbb{C}}

\newcommand{\X}{R_{n,G}}
\newcommand{\Y}[1]{R_{n,SU(#1)}}

\title{Cohomology of the space of commuting $n$-tuples in a compact Lie group}
\author{Thomas John Baird}
\date{}
\begin{document}
\maketitle

\begin{abstract}
Consider the space $Hom(\Z^n, G)$ of pairwise commuting $n$-tuples of
elements in a compact Lie group $G$. This forms a real algebraic
variety, which is generally singular. In this paper, we construct a
desingularization of the generic component of $Hom(\Z^n,G)$, which
allows us to derive formulas for its ordinary and equivariant
cohomology in terms of the Lie algebra of a maximal torus in $G$ and
the action of the Weyl group.  This is an application of a general
theorem concerning $G$-spaces for which every element is fixed by a
maximal torus.
\end{abstract}

\section{Introduction}

Moduli spaces of flat bundles are important in physics, where they
form critical level sets of Lagrangians for a number of important
quantum field theories, such as the Yang-Mills and Chern-Simons
theories. In this paper, we consider the moduli space of flat $G$-bundles
over a compact torus, $(S^1)^n$, where $G$ is a compact Lie group.  The fundamental group of $(S^1)^n$ is isomorphic to $\Z^n$,
and the holonomy map allows us to identify the moduli space of flat
bundles with $Hom( \Z^n, G)/G$, the space of homomorphisms of the
fundamental group into $G$, modulo conjugation. When n=2, $Hom( \Z^2, G)/G$ is isomorphic to the moduli space of semistable holomorphic $G_{\C}$
bundles over an elliptic curve, where $G_{\C}$ denotes the
complexification of $G$. The $G$-space $Hom( \Z^n, G)$ is the moduli space of based flat
connections and will be the principal object of study in this paper.

Because a homomorphism is determined by where it sends generators, we
may also identify $Hom( \Z^n, G)$ with $\{ (g_1,...,g_n)\in G^n | g_i g_j = g_j
g_i \ \forall i,j\}$, the space of pairwise commuting $n$-tuples in $G$.
$G$ carries a unique structure as a real algebraic group and $Hom( \Z^n, G)$
inherits the structure of a compact affine real algebraic variety. $G$
acts regularly on $Hom( \Z^n, G)$ by conjugation, and the quotient $Hom( \Z^n, G)/G$ is a space with
orbifold singularities.

A mistaken assumption about the nature of these spaces led to an
undercount of the number of vacuum states in supersymmetric Yang-Mills
theory over spatial $(S^1)^3$ by Witten in \cite{wi1}.  The source of
the error was identified in \cite{wi2}, and motivated a systematic
study of the spaces $Hom( \Z^n, G)/G$ in \cite{bfm} and \cite{ks}, with particular
focus on the case $n=3$. In both these papers, $Hom( \Z^n, G)/G$ is shown to be
describable in terms of the combinatorics of root systems. One
striking result to emerge is that $Hom( \Z^n, G)/G$ need not be connected, even
when $G$ is simply connected.

The study of the space $Hom( \Z^n, G)$ was initiated more recently in \cite{ac} in the broader context of Lie subgroups of $Gl_n(\C)$.
Their study was motivated by connections to orbifolds and pure braid
groups.  Their method is to consider the filtration of $Hom( \Z^n, G)$ by
subspaces $S_n(j,G) := \{(g_1,...,g_n)| $at least $j$ entries equal
$1_G \}$, where $1_G$ denotes the identity in $G$.  They show that
after suspending once, $Hom( \Z^n, G)$ decomposes up to homotopy equivalence as a
wedge of spaces:

\begin{equation}
\Sigma( Hom( \Z^n, G)) \sim \bigvee_{1\leq k \leq n} \Sigma( \bigvee_{{n
\choose k}} Hom( \Z^k, G)/S_k(1,G))
\end{equation}
They proceed to compute homology groups explicitly for $Hom( \Z^2, SU(2))$,
and $Hom( \Z^3, SU(2))$, though this method has so far not led to general formulas for cohomology.  See \cite{tgs} for a study of the case $G = SO(3)$,
including a description of connected components, fundamental groups
and $\Z_2$ homology.

Let $\X$ denote the connected component of
$Hom( \Z^n, G)$ containing the $n$-tuple $(1_G,...,1_G)$. $\X$ can be
characterized as those $n$-tuples whose entries lie in a common
maximal torus. If $G$ is connected and $g \in G$ is generic (i.e. has
centralizer a maximal torus), then any commuting $n$-tuple containing
$g$ lies in $\X$, so we will call $\X$ the generic component of $Hom(\Z^n,G)$. In many interesting cases (see Theorem \ref{class}), $Hom( \Z^n, G)$ is connected and so $X = Hom( \Z^n, G)$.

As motivation to focus on this component, we remark that the space $R_{2n,G}$ is isomorphic to the moduli space of based flat $G$-bundles over a genus $n$ Riemann surface, whose holonomy is reducible to a maximal torus, and that  $\X$ is also closely related to the corresponding space for the closed nonorientable surface $\R P^{\# (n+1)}$. In upcoming work (\cite{b}), the results of this paper will be used in the computation of the cohomology of moduli spaces of flat $SU(2)$-bundles over nonorientable surfaces.

In this paper, we construct a desingularization of $\X$, and present a nice
formula for its cohomology ring over fields with characteristic
relatively prime to the order of the Weyl group of $G$.

Our formula is a generalization of an old formula (see \cite{ghv}) for the cohomology
of $G = Hom(\Z,G)$,

\begin{equation}
H(G) \cong H(G/T \times T)^W
\end{equation}
where $T$ is a
maximal torus in $G$ and $W = N(T)/T$ is the Weyl group. This formula
is a consequence of the observation that the map $\phi: G/T \times T
\rightarrow G$, which sends $([g],t)$ to $\phi( ([g],t)) = gtg^{-1}$
is ``almost'' a covering map with covering transformation group $W$
in a sense we make precise in Section 2.  Basically, this means that $W$ acts freely on $G/T \times T$ leaving $\phi$ invariant and that $\phi^{-1}(g)/W$ is cohomologically equivalent to a point for all $g \in G$.

For example, when $G =SU(2)$, $G/T \times T$ is isomorphic to $S^2
\times S^1$ and $\phi^{-1}(g)$ is two points unless $g = \pm 1_G$ when
$\phi^{-1}(g) \cong S^2$. The quotient by $W = \Z_2$ is $\R P^2$ over
these exceptional points, which is cohomologically trivial for coefficient fields of characteristic not equal to 2.  Indeed, the space $S^2 \times_{\Z_2}S^1$ is
just the real blow up of $SU(2)$ over the points $\pm 1_G$.

For general $n$, the formula is

\begin{equation}
H(\X) \cong H(G/T \times T^n)^W
\end{equation}
with the same restriction on the coefficient field.  The
method of proof is to generalize this construction of an ``almost''
covering map to wider class of $G$-spaces that includes $\X$. The
main condition for this construction to work is that every stabilizer
must contain a maximal torus of $G$. As a bonus we also
get a formula for the equivariant cohomology,

\begin{equation}
H_G(X_{n,G}) \cong H_T(T^n)^W.
\end{equation}

In the case $G = SU(2)$ and $n \geq 2$, $\X$ has $2^n$ isolated
singularities at the points $( \pm 1_G )^n$.  Locally these
singularities are cones over $S^2 \times_{\Z_2}S^{n-1}$, where $\Z_2$
acts antipodally on each factor, so these singularities look like the
total space of the vector bundle $S^2 \times_{\Z_2} \R^n \rightarrow
S^2/\Z_2 = \R P^2$, quotiented by the zero section. $ G/T \times T^n
\cong S^2 \times (S^1)^n$, and its quotient $ S^2 \times_{\Z_2} (S^1)^n$
is a smooth manifold and can be identified as the blow up of $X$ by
these zero sections.


Some of the advantages of our approach over the approach in \cite{ac} is that we do obtain general formulas for the cohomology of $X$, and these formulas actually describe $H(X)$ as a ring, not just as a group. One of the disadvantages is that while we do get some negative results about torsion ($H(X)$ has no $p$-torsion for primes $p$ that do not divide the order of the Weyl group), the remaining torsion is not accessible by these methods.

The layout of this paper is as follows:

In Section 2, we define the notion of a cohomological principal
bundle.  When the structure group is finite, these are the ``almost"
covering maps described above, which have the important property that
the cohomology ring of the base is isomorphic to the invariant subring
of the cohomology of the total space.  I got the idea for this from
\cite{clm} where a construction similar to the one in this paper is used in the study of the
moduli space of flat $SU(2)$ connections on a Riemann surface.
Cohomological principal bundles also emerge in a different way in the study of $SU(2)$
connections on nonorientable surfaces in \cite{b} and seem to be a
recurring phenomenon in such problems.

In Section 3, we prove a general result (Theorem \ref{thm1}) concerning the cohomology of spaces admitting $G$ actions
for which every point is fixed by a maximal torus.

In Section 4, we apply the theorems proved in Section 3 to the case of
commuting $n$-tuples, obtaining formulas for the ordinary and
equivariant cohomology of $\X$.

In Section 5, we explore some consequences of these formulas by
working out a more explicit description in the case $G=SU(2)$ and computing Poincar\'e
polynomials for $G = SU(2)$, $SU(3)$ and $SU(4)$.

Also included are two appendices reviewing some results that are
applied in the paper. Appendix A is an overview of equivariant
cohomology from the Borel model perspective, and in Appendix B we work
out formulas for $H(T^n)$ and $H(G/T)$ as modules of the Weyl group.

I would like to thank Lisa Jeffrey and Paul Selick for all of their
suggestions and guidance, as well as Fred Cohen and Eckhard Meinrenken
for helpful conversations.

\section{Cohomological Principal Bundles}

Let $ f: X \rightarrow Y $ be a continuous map between topological
spaces $X$ and $Y$, and let $ \Gamma $ be a topological group
acting freely on the right of $X$, such that $ X \rightarrow X /
\Gamma $ is a principal bundle.

\begin{define}\label{princ}
We say $(f:X  \rightarrow Y, \Gamma )$ is a cohomological
principal bundle for the cohomology theory H if:
\\
i) $f$ is a closed surjection
\\
ii) $f$ descends through the quotient to a map $h$,
\\
\begin{equation}\begin{CD}
\xymatrix{ X \ar[d]^{\pi} \ar[dr]^f \\
       X/\Gamma \ar[r]^h & Y}
\end{CD}\end{equation}
\\
iii) $H(h^{-1}(y)) \cong H(pt)$ for all $y \in Y$
\end{define}

Let $H(X,F)$ denote sheaf cohomology of the constant sheaf $F_X$, where $F$ is a field (in all applications we have in mind, sheaf cohomology is isomorphic to singular cohomology). To prove Proposition \ref{invt} we will require the following two standard results (see \cite{br} 11.7 and 19.2):

\begin{thm}\label{VB}[Vietoris-Begle mapping theorem] Let $h: Z \rightarrow Y$ be a closed surjection, where $Z$ is a paracompact Hausdorff space.  Suppose that for all $y \in Y$, $H(f^{-1}(y),F) \cong H(pt,F)$. Then

\begin{equation}
f^*: H(Y,F) \rightarrow H(Z,F)
\end{equation}
is an isomorphism.
\end{thm}

\begin{thm}\label{cov} Let $X$ be a topological space, let $\Gamma$ be a finite group acting on $X$ and let $\pi: X \rightarrow X/\Gamma$ denote the quotient map onto the orbit space $X/\Gamma$.  If $F$ is a field satisfying $gcd(char(F),\#\Gamma)=1$, then

\begin{equation}
\pi^*: H(X/\Gamma,F) \rightarrow H(X,F)^{\Gamma}
\end{equation}
is an isomorphism, where $H(X,F)^{\Gamma}$ denotes the ring of $\Gamma$ invariants.
\end{thm}


Combining Theorem \ref{VB} and Theorem \ref{cov}, we get:

\begin{prop}\label{invt}  Let $\Gamma$ be a finite group of order N, let $X$ be a paracompact Hausdorff space and
let $(f:X  \rightarrow Y, \Gamma )$ be a cohomological principal
bundle for $H(.,F)$, where
gcd(char(F),N) = 1. Then $f^*$ induces an isomorphism

\begin{equation}
f^*: H(Y,F) \cong H(X,F)^{\Gamma}
\end{equation}
where $H(X,F)^{\Gamma}$ denotes the ring of $\Gamma$ invariants.
\end{prop}

\begin{proof}
Since $X$ is a paracompact Hausdorff space, $X/\Gamma$ is as well.  The induced map $h: X/\Gamma \rightarrow Y$ satisfies the conditions of the Vietoris mapping theorem, so

\begin{equation}
h^*: H(Y,F) \cong H(X/\Gamma,F)
\end{equation}

On the other hand, by Theorem \ref{cov}, $\pi^*: H(X/\Gamma,F) \cong H(X,F)^{\Gamma}$.  Since $f^* = \pi^* \circ h^*$ this completes the proof.
\end{proof}

Though Proposition \ref{invt} suffices for our purposes, it will more convenient to use the following corollary.

\begin{cor}\label{corollary}Let $\Gamma$ be a compact Lie group with N connected components,
and let $\Gamma_0$ be the identity component.  Let $(f:X \rightarrow Y, \Gamma )$ be a cohomological principal bundle for $H(.,F)$, where $gcd(char(F),N) = 1$ and $X$ is a paracompact Hausdorff space. Then
$H(Y,F) \cong H(X/\Gamma_0,F)^{\Gamma/\Gamma_0}$.
\end{cor}
\begin{proof}
$f$ descends to a map $g: X/\Gamma_0 \rightarrow Y$ and the
residual action of $\Gamma / \Gamma_0$ acts on $X/ \Gamma_0$
making the pair into a cohomological covering map.  Since $X \rightarrow X/\Gamma_0$ is a principal bundle, it is a closed map and thus $X/\Gamma_0$ is paracompact Hausdorff (see \cite{e} 5.1). The result then follows from Proposition \ref{invt}.
\end{proof}

\section{Cohomology of G-spaces with stabilizers containing a maximal
torus}

Let $G$ be a connected compact group, $T$ a maximal torus of $G$ and $X$ a space on which $G$ acts.
If every point $x \in X$ is fixed under this action by some maximal
torus of $G$, then because all maximal tori in $G$ are conjugate, every $G$ orbit must intersect the $T$ fixed point set $X^T$. It follows that the map

\begin{equation}
\phi : G \times X^T \rightarrow X,  \phi((g,x)) = g \cdot x
\end{equation}
is surjective.  $G$ acts on $G \times X^T$ by $g \cdot (h,x) =
(gh,x)$, and $\phi$ is equivariant for this action. The normalizer of
$T$ in $G$, denoted $N(T) = N_G(T)$, acts freely on $G \times X^T$ from the right by

\begin{equation}
(g,x) \cdot n = (gn, n^{-1} \cdot x)
\end{equation}
leaving $\phi$ invariant and commuting with the $G$ action. We will
show in Theorem \ref{thm1} that under very mild conditions
the pair $( \phi: G \times X^T \rightarrow X, N(T))$ is a
cohomological principal bundle. We begin with a couple of lemmas.

\begin{lem}\label{lem1}
Let $G$ act on $X$ from the left and let $x \in X^T$.  Then $g \cdot x
\in X^T$ if and only if $g \in N(T)G_x^0$, where $G_x^0$ is the
identity component of the stabilizer $G_x$.
\end{lem}

\begin{proof}
If $g \cdot x \in X^T$, then $g^{-1} t g \cdot x = x$ for all $t \in
T$, so

\begin{equation}
 g^{-1}Tg \subset G_x
\end{equation}
Since $T$ is maximal in $G$, it is also maximal in $G_x$, so for
some $h \in G_x^0$, $h^{-1}g^{-1}Tgh = T$, and thus $g \in N(T)G_x^0$.
The other direction is clear.
\end{proof}

Let $W_G := N_G(T)/T $ denote the Weyl group of $G$.

\begin{lem}\label{lem2}
Let $( \phi: G \times X^T \rightarrow X, N(T))$ be defined as above.
For every $x \in X$, $H(\phi^{-1}(x)/N(T),F) \cong H(pt,F)$, for $F$ satisfying $gcd(char(F),\# W_G)=1$.
\end{lem}

\begin{proof}
We may assume by equivariance that $x \in X^T$.

\begin{equation}
\phi^{-1}(x) = \{ (g,y) | g \cdot y = x \} = \{(g,y) | y \in X^T,
g^{-1} x = y \} \cong G_x^0N(T)
\end{equation}
where this last isomorphism follows from the preceding lemma.  It
follows that

\begin{equation}
\phi^{-1}(x)/N(T) \cong G^0_xN(T)/N(T) \cong G_x^0/N_{G_x^0}(T)
\end{equation}
Now since $\#W_{G_x^0}$ divides $\#W_G$, we deduce from Proposition \ref{hsi} that

\begin{equation}
H(G_x^0/N_{G_x^0}(T),F) \cong H(pt,F)
\end{equation}
completing the proof.
\end{proof}

\begin{thm}\label{thm1}
Let $G$ be a compact, connected Lie group with maximal torus $T$,
acting on a paracompact Hausdorff space $X$.  Suppose that
for every point $x \in X$, $G_x$ contains a maximal torus of $G$.
Then $( \phi: G \times X^T \rightarrow X, N(T))$ is a cohomological principal bundle for H(.,F), where $gcd(char(F),\# W_G) = 1$.  In particular:

\begin{equation}
H(X,F) \cong H(G/T \times X^T,F)^{W_G}
\end{equation}

\end{thm}

\begin{proof}
First note that since $X^T$ is a closed subset of $X$ it inherits a paracompact Hausdorff topology.  The only conditions in Definition \ref{princ} that are not immediate are (iii), which follows from Lemma \ref{lem2}, and closedness of the map $\phi$. But $\phi$ is a restriction of the action map $G \times X \rightarrow X$ which is easily shown to be closed, so $\phi$ is also closed.

The final assertion follows from Corollary \ref{corollary}

\end{proof}

It is helpful to consider two extreme cases of this theorem:

\begin{example}
Suppose that $G$ acts trivially on $X$.  Then $X = X^T$ and $G \times_{N(T)}X^T = G/N(T) \times X$. By Proposition \ref{hsi}, we know that $G/N(T)$ has trivial cohomology, so by the Kunneth theorem $H(X,F) \cong H(G/N(T) \times X,F) \cong H(G/T \times X^T,F)^{W_G}$.
\end{example}

\begin{example}
Suppose that $G$ acts so that $G_x$ is a maximal torus for every $x\in X$. By Lemma \ref{lem1} we deduce that each orbit must intersect $X^T$ precisely $\#W_G$ times. It follows easily that $G/T \times X^T \rightarrow X$ is a covering space map with deck transformation group $W_G$, and thus $H(G/T \times X^T,F)^{W_G} \cong H(X,F)$.
\end{example}

\begin{cor}\label{haha}
Let $G$ act on $X$ satisfying the hypotheses of Theorem \ref{thm1}.  Then $dim H(X,\C) = dim H(X^T,\C)$.
\end{cor}

\begin{proof}
It is a general property proven using characters, that if $V$ is a
finite dimensional representation of a finite group $\Gamma$, then
$dim ( \C \Gamma \otimes V)^{\Gamma} = dim V$. Because $H(G/T,\C) \cong \C W$ as a left W representation (see Proposition \ref{big}) it follows that:

\begin{equation}
dim H(X,\C) = dim (H(G/T,\C) \otimes H(X^T,\C))^W = dim H(X^T,\C)
\end{equation}

\end{proof}

When $X$ is a smooth compact manifold or more generally when the set of infinitesimal stabilizers of the action is finite, Corollary \ref{haha} combined
with Proposition \ref{formal2} implies that the action is $T$
equivariantly formal (see Appendix A for a review of equivariant cohomology).  We adopt the convention that equivariant cohomology is assumed to be taken with complex coefficients unless otherwise stated, i.e. $H_G(.) = H_G(.,\C)$.

Now if we have a $G$ equivariant map $\rho: Y \rightarrow X$ between
$G$-spaces $Y$ and $X$, we obtain a map between homotopy quotients
$\rho_G:Y_G \rightarrow X_G$.  If $\rho$ induces an isomorphism $H(X,\C) \cong
H(Y,\C)$, we see by considering the Serre spectral sequences of the
standard fibrations of $X_G$ and $Y_G$ over $BG$ that $\rho_G^*:
H_G(X) \cong H_G(Y)$.

\begin{thm}\label{core}
Let $G$ act on $X$ satisfying the hypotheses of Theorem \ref{thm1}.  Then $H_G(X) \cong H_T(X^T)^W = (H(X^T)
\otimes H_T)^W$.
\end{thm}

\begin{proof}
The $G$-equivariant map $\phi : G/T \times_W X^T \rightarrow X$ induces an isomorphism
in $\C$ cohomology, so

\begin{equation}
H_G(X) \cong H_G(G/T \times_W X^T) \cong H_G(G/T \times X^T)^W
\end{equation}
It follows from a well known formula in equivariant cohomology (equation (\ref{stan})) that $H_G(G/T \times X^T) \cong H_T(X^T)$, so
we deduce that

\begin{equation}
H_G(X) \cong H_T(X^T)^W
\end{equation}
\end{proof}

\begin{rmk}
This theorem can also be proven by showing directly that $(EG \times
X^T \rightarrow X_G, N(T))$ is a cohomological principal bundle.
\end{rmk}

\begin{rmk}
By the familiar identity $H_G(X) \cong H_T(X)^W$, Theorem \ref{core} shows that the localization map $i^*: H_T(X) \rightarrow H_T(X^T)$ restricts to an isomorphism between the Weyl invariant subrings.
\end{rmk}

Just to clarify this theorem, we work out carefully how $W$ acts on
$H_T(X^T)$.  If $E$ is the total space of a universal $G$-bundle, then
$(G/T \times X^T)_G = E \times_G (G/T \times X^T)$.  If $n \in N(T)$
represents an element of $W$, and $(e,g,x)$ represents an element of
$(G/T \times X^T)_G$, then $ [(e,g,x)]\cdot [n] = [(e,gn,n^{-1} \cdot
x)]$.  Thus in terms of the identification with $E \times_T X^T$ the
action looks like $[(e,x)] \cdot [n] = [(e \cdot n, n^{-1} \cdot x)]$.
If we turn the right representation of $W$ on $H(X^T)$ into a
left representation in the usual way, the action of $W$ on $H_T(X^T)$
is just the tensor product of the representations on $H(X^T)$ and
$H_T$.

\section{Cohomology of $\X$}

Let $G$ be a connected, compact Lie group and let $\X$ be the identity component of $
\{(g_1,...,g_n) \in G^n| g_ig_j = g_jg_i \ \forall i,j \} \cong
Hom(\Z^n,G)$, topologized as a subspace of $G^n$.  The following theorem, quoted from \cite{ks}, shows that in many cases $\X$ is in fact the only component of $Hom(\Z^n,G)$

\begin{thm}\label{class}
Let $G$ be a compact, simple Lie group. The space $Hom(\Z^n,G)$ is connected if and only if any of the following conditions are met:

i) $n=1$, and $G$ is connected.

ii) $n=2$ and $G$ is 1-connected.

iii) $n \geq 3$ and $G = SU(m)$ or $Sp(m)$, $m \geq 1$.
\end{thm}

The following lemma, which is a consequence of classification of components of $Hom(\Z^n,G)$ in \cite{bfm} and \cite{ks}, gives a concrete description of $\X$.

\begin{lem}\label{torus}
For a compact Lie group G, every commuting n-tuple in $\X$ lies in a maximal torus of G.
\end{lem}



It follows that every point in $\X$ is fixed by a maximal torus under the conjugation action.

Our main result is:

\begin{thm}\label{consult}
Let G be a connected, compact Lie group and let $T$ be a maximal torus
in G. The pair ($\phi :G \times T^n \rightarrow \X, N(T))$ forms a
cohomological principal bundle for cohomology over fields $F$ of
characteristic relatively prime to $\# W$, where $W$ is the Weyl group.
In particular, $H(\X,F) \cong H(G/T \times T^n,F)^W$.
\end{thm}

\begin{proof}

This is a straightforward application of Theorem \ref{thm1}.







\end{proof}

We deduce using Theorem \ref{core}:

\begin{cor}\label{aaaaa}
$H_G(\X) \cong H_T(T^n)^W$.
\end{cor}

\begin{rmk}
It may be shown (see \cite{s}), that the quotient $G \times_{N(T)}T^n$ is a nonsingular real algebraic variety.  The induced map
$h: G \times_{N(T)}T^n  \rightarrow \X$ is surjective,
regular and induces birational equivalence.  Thus $h: G \times_{N(T)}T^n
 \rightarrow \X$ is a resolution
of singularities for $\X$.
\end{rmk}

\begin{rmk}
$G$ acts on $G \times_{N(T)}T^n=G/T \times_W T^n$ by left
multiplication on the $G$ factor, making $h$ into a $G$-equivariant
map.  $h$ descends to an isomorphism between $G \backslash (G/T
\times_W T^n) = T^n/W$ and $\X/G$.
\end{rmk}

Since $G/T$ has a CW-structure with only even dimensional Schubert
cells, $H(G/T, \Z)$ is torsion free, and $H(T^n, \Z)$ is also
torsion free. If $(char(F),\#W) = 1$, it follows by Lemma \ref{annoy}
that $dim H^k(\X,F) = dim H^k(\X, \C)$. The universal coefficient
theorem now applies to prove: 

\begin{cor}\label{conf}
$Tor(\Z_p,H^*(\X, \Z)) =0$, for primes $p$ satisfying $(p, \#W)=1$.
\end{cor}

\begin{lem}\label{annoy}
Let $\Gamma$ be a finite group acting linearly on a free, finitely generated
$\Z$-module $M \cong \Z^n$.  Given any field $F$, there is an
induced linear action on $M_F = M\otimes F$. If $gcd( char(F), \#\Gamma)=1$, then $dim_F( (M_F)^{\Gamma}) = rank_{\Z}(M^{\Gamma})$
\end{lem}

\begin{proof}
Because $M^{\Gamma}$ is a saturated sublattice of $M$(in the sense that $\alpha m \in M^{\Gamma}$ for nonzero $\alpha \in \Z$ implies $m \in M$),  it follows that any basis of $M^{\Gamma}$ extends to a basis of $M$.

Let $\psi: M \rightarrow M\otimes F$ sends $ m$ to $m \otimes 1$ and let
$V := span_F\{ \phi(M^{\Gamma})\}$. Then $dim_F(V) = rank_{\Z} M^{\Gamma}$ and $V \subset (M_F)^{\Gamma}$. In fact, we will show $V = (M_F)^{\Gamma}$, thus proving the theorem.

Because $\#\Gamma$ is invertible in $F$, we can define a projection operator $P: M_{F} \rightarrow M_F^{\Gamma}$ by

\begin{equation}
P(x) = \frac {1} {\# \Gamma} \sum_{g\in \Gamma}gx
\end{equation}
The image of $P$ is spanned by vectors $\sum_{g\in \Gamma}g \psi(m)= \psi(\sum_{g\in \Gamma}gm)$, which all lie in $V$ which proves that $V=(M_F)^{\Gamma}$ completing the proof.

\end{proof}

In view of Corollary \ref{conf}, little will be lost by focusing on the characteristic zero case from now on.

The first homology group $H_1(T^n,\R)$ is canonically isomorphic to the Lie algebra, $Lie(T^n) = \lie{t^n}$.  It follows that $H^*(T^n) \cong \bigwedge \lie{t}^{*n}$ with its usual grading.  It is also well known that the equivariant cohomology of a point $H_T(pt, \R)$ is canonically isomorphic to the symmetric algebra $S \lie{t}^*$, with grading $deg( S^i \lie{t}) = 2i$ (see Appendix A for more details).  Combining this with Theorem \ref{consult} and Corollary \ref{aaaaa} results in the following simple formulas:

\begin{equation}
H_G(\X, \R) \cong (\bigwedge \lie{t}^{*n} \otimes S\lie{t}^*)^W
\end{equation}
and
\begin{equation}
H(\X,\R) \cong (\bigwedge \lie{t}^{*n}
\otimes S\lie{t}^*)^W / <S\lie{t}^{*W+}>
\end{equation}
where $<S\lie{t}^{*W+}>$ denotes the ideal generated by the image of
the ring of positive degree elements in $S\lie{t}^{*W} \cong H_G$.
\section{Examples}

Throughout this section we will always assume cohomology to be taken
with coefficients in $\C$.

\subsection{$G = SU(2)$}

We begin with the case $G =SU(2)$.  In this case both $H(\X)$ and
$H_G(\X)$ admit simple, explicit descriptions.  Here, $\lie{t}$ is one
dimensional so that $\lie{t}^{*n} \otimes \C \cong \C^n$ and $W \cong \Z_2$ acts
by multiplication by $\{ \pm 1\}$.  $G/T \cong \C P^1$ and $W$ acts as
an orientation reversing involution. Letting $y$ denote a generator of
$H^2(G/T)$, we have:

\begin{equation}
H(G/T \times T^n) \cong \bigwedge \C^n \otimes \frac{\C [y]}{ <y^2
=0>}
\end{equation}
where $y$ has degree 2 and $\C^n$ has degree 1.  Identifying
$H(X_{n,SU(2)})$ with the Weyl invariant subring gives

\begin{equation}
H(X_{n,SU(2)}) \cong (\oplus_{i\  even}\mathord{\bigwedge}^i \C^n) \oplus
(\oplus_{i ~odd} \bigwedge^i \C^n \otimes y)
\end{equation}


We can easily compute Betti numbers to be $dim H^d(X_{n,SU(2)}) = {n
\choose d}$ if $d$ is even, and ${n \choose d-2}$ if $d$ is odd.

Similarly, the $SU(2)$ equivariant cohomology of $X_{n,SU(2)}$ is
isomorphic to

\begin{equation}
H_G(X_{n,SU(2)}) \cong \oplus_{i+j \equiv 0 ~mod 2} (\bigwedge^i\C^n)
\otimes y^j
 \subseteq \bigwedge\C^n \otimes \C[y]
\end{equation}
where the generators of the exterior algebra have order 1 and that of
the polynomial algebra has order 2.

\subsection{Poincar\'e polynomials}

The computation of the Poincar\'e polynomial of $\X$ can be organized
using polynomials with character coefficients.  If a group $W$ acts on
a graded vector space $\oplus V_i$ with character $\chi^i$ on $V_i$,
then the Poincar\'e polynomial of the representation is $\sum_i \chi^i
t^i$.  Direct sums and tensor products of graded representations pass to
sums and products of their Poincar\'e polynomials.  The ordinary
Poincar\'e polynomial of the invariant part of the representation can
be extracted by taking the inner product of the coefficients with the
character of the trivial representation.

For example, when $G = SU(2)$, the character table for $W = S_2$ is:

\begin{tabular}{c|ccc}
 & (1) & (12) \\
\hline
$\chi_1$ & 1 & 1 \\
&&&\\
$\chi_2$ & 1 & -1 \\
\end{tabular}
\\
\\
and the formula for the Poincar\'e polynomial of $\Y{2}$ is

\begin{align*}
&P_t(H(\Y{2})) = <\chi_1, (\chi_1 + \chi_2t)^n(\chi_1 +
\chi_2t^2)>\\
&= 1/2 [ (1+t)^n(1+t^2) + (1-t)^n(1-t^2)] =1/2 \sum_i {n \choose
i}[(1+(-1)^i) t^i + (1-(-1)^i)t^{i+2}]
\end{align*}
which concurs with our earlier calculation.

For $G=SU(3)$, the character table for $W = S_3$ is

\begin{tabular}{c|ccc}
  & (1) & (12)&(123)\\
\hline
$\chi_1$  & 1&1&1 \\
&&&\\
$\chi_2$ &1&-1&1 \\
&&&\\
$\chi_3$ &2&0&-1 \\
\end{tabular}
\\
\\
and the Poincar\'e polynomial is

\begin{align*}
&P_t(H(\Y{3})) = < \chi_1, (\chi_1 + \chi_3 t + \chi_2
t^2)^n(\chi_1 + \chi_3 t^2 + \chi_3 t^4 + \chi_2 t^6)> \\
&= 1/6[ (1+t+t^2)^n(1+t^2+t^4+t^6) + 3(1-t^2)^n(1-t^6) +
2(1-t+t^2)^n(1-t^2-t^4+t^6)]
\end{align*}

For $G = SU(4)$, the character table for $W = S_4$ is

\begin{tabular}{c|ccccc}
  & (1) & (12)&(123)&(1234)&(12)(34)\\
\hline
$\chi_1$  & 1&1&1&1&1 \\
&&&\\
$\chi_2$ &1&-1&1&-1&1 \\
&&&\\
$\chi_3$ &2&0&-1&0&2 \\
&&&\\
$\chi_4$ &3 &1& 0& -1& -1 \\
&&&\\
$\chi_5$ & 3 & -1& 0& 1&-1 \\
\end{tabular}
\\
\\
and the Poincar\'e polynomial of $\Y{4}$ is:

\begin{align*}
 &P_t(H(\Y{4})) = < \chi_2, (\chi_1 + \chi_4t + \chi_5t^2
+\chi_2t^3)^n(\chi_1 + \chi_4 t^2 + \\
&(\chi_3+\chi_4) t^4 + (\chi_4 +\chi_5) t^6 + (\chi_3+\chi_5)t^8 +
\chi_5 t^{10} + \chi_2 t^{12})>
\end{align*}

\appendix
\section{Equivariant Cohomology}

In this section, we briefly review equivariant cohomology from the
Borel model perspective. Our main sources are \cite{ab} and \cite{h}.

Let $G$ be a compact Lie group and $X$ a topological space with
continuous left $G$ action. We define the equivariant cohomology
$H_G(X)$ to be the singular cohomology over $\C$ of the space
$X_G$ obtained from a universal $G$-space $EG$ by the mixing
construction

\begin{equation}\label{mix}
X_G := EG \times_G X
\end{equation}

In equation (\ref{mix}), $G$ acts on $EG$ from the right, and $EG
\times_G X$ denotes the quotient of $EG \times X$ by the relation
$(e,gx) \sim (eg,x)$.  $X_G$ projects naturally onto the
classifying space $BG = EG/G$ forming a fibre bundle with fibre
$X$. We will call this the standard fibre bundle.

We adopt the notation $H_G := H(BG,\C)$. The standard fibre bundle
$\pi : X_G \rightarrow BG$ induces a natural map $ \pi^* : H_G
\rightarrow H_G(X)$, making $H_G(X)$ into an $H_G$-module.

The Serre spectral sequence of the standard fibration is frequently
used to compute $H_G(X)$.  If this spectral sequence collapses at the
$E_2$ page, then $H_G(X) \cong H(X) \otimes_{\C} H_G$ as $H_G$-modules
and we say that the action is formal.  When the action is formal, the
Leray-Hirsch theorem implies that

\begin{equation}\label{lerayhirsch} H_G(X)/<H_G^+> \cong H(X) \end{equation} where
$<H_G^+>$ denotes the ideal in $H_G(X)$ generated by the image of the
positive degree ideal in $H_G$ under $\pi^*$. We will use the
following simple criterion for formality.

\begin{prop}\label{formal}
Let $G$ be a compact, connected Lie group, and $X$ a space such that
$H^i(X)=0$ for i odd.  Then any action of $G$ on $X$ is formal.
\end{prop}
\begin{proof}
According to \cite{ab}, $H^i_G = 0$ for odd i, so that the Serre spectral
sequence of the standard fibration satisfies $E_2^{p,q} = 0$ unless
$p$ and $q$ are even. Thus all subsequent differentials must be zero
and the spectral sequence collapses at $E_2$.
\end{proof}

Let $T = (S^1)^m$ be a torus. Let $X$ be a smooth compact manifold on
which $T$ acts smoothly, or more generally let the set of infinitesimal stablizers of the action be finite.  The fixed point set $X^T$ includes into $X$
by the map $i$, inducing a map $i^*: H_T(X) \rightarrow H_T(X^T)$. The
localization theorem [AB] states:

\begin{thm}\label{loc}
The kernel and cokernel of $i^*$ are torsion $H_T$-modules.
\end{thm}

Combining the localization theorem with the Serre spectral sequence of
the standard fibration gives (see \cite{gs}):

\begin{prop}\label{formal2}
With conditions as in \ref{loc}, the action of T on X is formal if and
only if $dim H(X^T) = dim H(X)$.
\end{prop}

When $G$ is a connected, compact Lie group with maximal torus $T$,
any $G$-space $X$ becomes a $T$-space by restricting the action.
$EG$ also forms a model for $ET$ under the restricted action, and
there is a natural map $\rho: X_T \rightarrow X_G$, which is a
fibre bundle with fibre $G/T$.  In fact, there is a natural
identification

\begin{equation}\label{stan}
X_T \cong (X \times G/T)_G
\end{equation}
where $G$ acts by left multiplication on the homogeneous space
$G/T$, and $G$ acts on the product $X \times G/T$ via the diagonal
action. Under this identification, the fibre bundle is just the
projection map $(X \times G/T)_G \rightarrow X_G$.

The homotopy quotient $X_T$ retains a natural action of the Weyl
group $W = N(T)/T$ and we obtain an isomorphism:

\begin{equation}\label{yeesh}
\rho^*: H_G(X) \cong H_T(X)^W
\end{equation}

In fact by Proposition \ref{hsi}, $(\rho: X_T \rightarrow X_G, W)$ is a cohomological
principal bundle, so this formula follows from Proposition
\ref{invt}.

\begin{prop}\label{hsi}
Let G be a compact, connected Lie group with maximal torus T.  Then $H(G/N_G(T),F) \cong H(pt,F)$ for $gcd(char(F),\#W_G)=1$.
\end{prop}

\begin{proof}
Since $G/T$ has a CW structure with $\#W_G$ cells (the Schubert cells) all in even degrees,

\begin{equation}
H(G/N_G(T),F)=H(G/T,F)^{W_G}
\end{equation}
has Euler characteristic 1, and is zero in odd degrees, and so must be trivial.
\end{proof}

As a graded ring,

\begin{equation}\label{H_G}
H_G \cong (S \clie{t}^*)^W
\end{equation}
the Weyl invariant part of the symmetric algebra over
$\clie{t}^*$, where $\clie{t}$ is the complexified Lie algebra of
a maximal torus in $G$ and we set the degree of $\clie{t}^* = 2$.
This is a polynomial algebra with even degree generators. For instance, when
$G = U(N)$ the generators are the universal Chern classes
$\{ c_1, c_2,..., c_N \}$.

\section{$W$-module structure of $H(T^n)$ and $H(G/T)$}

Let $G$ be a compact, connected Lie group and $T$ be a maximal torus
in $G$. We assume throughout this section that $H(.)$ is singular
cohomology over $\C$.

Let $\lie{t}_{\C} = \lie{t} \otimes \C$ be the complexification of
the Lie algebra of $T$.  It is immediately verified by the reader
using invariant de Rham forms that

\begin{equation}\label{obvious}
H(T^n) \cong (\bigwedge \lie{t_{\C}^*})^{\otimes n}
\end{equation}
where the degree of $\clie{t}^*$ is 1. The Weyl group $W = N(T)/T$
acts on $T^n$ by conjugation, inducing a $W$-module structure on
$H(T^n)$ which is the same one induced by the usual action on
$\clie{t}^*$.

The W-module structure of $H(G/T)$ is somewhat more complicated. We will
content ourselves here with determining it as an ungraded
representation.

\begin{prop}\label{big}
As an (ungraded) left $W$-module, $H(G/T) \cong \C W$, the group
ring.
\end{prop}
\begin{proof}
The set $N(T)/T = W$ will come up in this proof both as a group,
and as a subset of $G/T$.  To avoid confusion, we use $W$ when we
regard it as a group and $N(T)/T$ when we regard it a subset.

$G$, and hence $T$ acts on $G/T$ by multiplication on the left. $G/T$ possesses a
CW structure of even dimensional Schubert cells, so by Proposition
\ref{formal} we know that the $T$ action is formal and that

\begin{equation}\label{eqn a}
H_T(G/T) \cong H(G/T) \otimes H_T
\end{equation}
as $H_T$-modules.  More precisely, the natural map $\phi : H_T(G/T)
\rightarrow H(G/T)$ is surjective. A linear section $s: H(G/T)
\rightarrow H_T(G/T)$ which respects the grading determines an
isomorphism $\tilde{s}: H(G/T) \otimes H_T \rightarrow H_T(G/T)$, by

\begin{equation}
\tilde{s} (x, \alpha) = s(x) \alpha
\end{equation}

The fixed point set of the $T$ action on $G/T$ is $N(T)/T$.
Because $H_T(G/T)$ is a free $H_T$-module, the localization map
$i^*: H_T(G/T) \rightarrow H_T(N(T)/T)$ is an inclusion whose
cokernel is a torsion $H_T$-module. Since $T$ acts trivially on
$N(T)/T$,

\begin{equation}\label{eqn b}
H_T(N(T)/T) \cong H(N(T)/T) \otimes H_T
\end{equation}
canonically as rings.  Now $W$ acts on the right of $G/T$,
preserving $N(T)/T$, and commuting with the $T$ action, so it
induces an action on $H_T(G/T)$ and $H_T(N(T)/T)$ which is
equivariant for $i^*$. By choosing the section $s$ to be
$W$-equivariant also, we find that in terms of the equations
(\ref{eqn a}) and (\ref{eqn b}), $W$ acts trivially on the $H_T$
factor. Let $Q_T$ denote the quotient field of $H_T$. Tensoring
equations (\ref{eqn a}) and (\ref{eqn b}) by $Q_T$ gives an
isomorphism of $W$ representations:

\begin{equation}
H(G/T) \otimes_{\C} Q_T \cong H(N(T)/T) \otimes_{\C} Q_T
\end{equation}

But this proves the theorem, because $Q_T$ is a field extension
of $\C$
\end{proof}

For an alternative proof, see \cite{cg} chapter 6.

\begin{thm}
As graded rings and W modules, $H(G/T) \cong H_T/<H_T^{W+}>$
\end{thm}

\begin{proof}
H(G/T) is trivial in odd degrees, so by equation (\ref{lerayhirsch})
and Proposition \ref{formal} we know
\begin{equation}
H(G/T)\cong H_G(G/T)/< H_G^+>
\end{equation}
On the other hand, by equation (\ref{stan}) we know $H_G(G/T) \cong
H_T$, and by equation (\ref{yeesh}) that $H_G^+ \cong H_T^{W+}$, which
verifies the formula.
\end{proof}

We point out a consequence of this theorem.

\begin{cor}
$S \clie{t}^* \cong \C W \otimes (S \clie{t}^*)^W$ as an ungraded
W representation, under the usual left action.
\end{cor}

\end{document}